\documentclass[10pt]{article}
\usepackage{amsmath,amssymb,amsthm,graphicx,hyperref,amsrefs}

\newcommand{\MRnumber}[1]{~\href{http://www.ams.org/mathscinet-getitem?mr=#1}{{\bf MR~#1}}}

\newcommand{\NN}{{\mathbb N}}
\newcommand{\ZZ}{{\mathbb Z}}

\newcommand{\floor}[1]{\left\lfloor #1 \right\rfloor}

\newcommand{\fp}[1]{\{ #1 \}}
\newcommand{\tf}[1]{[\![ #1 ]\!]}

\newcommand{\B}{B}
\newcommand{\BB}{B^\prime}

\newtheorem{lem}{Lemma}
\newtheorem{theorem}{Theorem}

\title{On the Sum of the Heights of Sturmian Factors}
\author{Kevin O'Bryant}

\begin{document}
\maketitle\thispagestyle{empty}

\begin{abstract}
A binary word is a map $W\colon \NN \to \{0,1\}$, and the set of factors of $W$ with length $n$ is
$F_n(W):=\{\big(W(i),W(i+1),\dots,W(i+n-1)\big)\colon i\ge 0\}$. A word is {\em Sturmian} if
$|F_n(W)|=n+1$ for every $n\ge1$. We show that the sum of the heights (also known as hamming
weights) of the $n+1$ factors with length $n$ of a binary Sturmian word has the same parity as $n$,
independent of $W$.
\end{abstract}

Many facts are known about the factors of length $n$ of a Sturmian word $W$. Among the many
noteworthy results are: that $F_n(W)$ is closed under reversals (the map that takes
$(w_1,\dots,w_n)$ to $(w_n,\dots,w_1)$) \cite{Lothaire}*{Prop 2.1.19}; that the volume of the
convex hull of $F_n(W)$ is $1/n!$, independent of $W$ \cite{OBryant}*{Thm 1.1}; and that as $W$
varies over all Sturmian words, $F_n(W)$ takes on $\sum_{i=1}^n \phi(i)$ values. We direct the
reader to either \cite{Shallit}*{Chap 9} or \cite{Lothaire}*{Chap 2} for an introduction to
Sturmian words. To these we add
    \begin{theorem} \label{thm:main}
    For every binary Sturmian word $W$ and every positive integer $n$,
    \[ \sum_{\vec w \in F_n(W)} h(\vec w) \equiv n \pmod2,\]
    where $h(\vec w)$ is the number of components of $\vec w$ that are `1'.
    \end{theorem}

A natural approach to proving this is to observe that since $F_n(W)$ is closed under reversal, we
can pair off non-palindrome factors that have the same height $h(\vec w)$, and therefore it
suffices to consider only the palindromes in $F_n(W)$. Moreover, if a palindrome has even length,
then it must have even height, and so the `even-$n$ case' of our theorem does follow easily from
the `closure under reversal' property. When $n$ is odd, the situation is more complicated as a
palindrome may have even or odd weight, and there are always two~\cite{Droubay}:
    \[\big\{ (1,0,1,0,1,0,1),(1,1,0,1,0,1,1)\big\} \subseteq F_7(c_{1/\sqrt3}),\]
where $c_{1/\sqrt{3}}$ is a particular Sturmian word defined below. Our proof does not follow this
line, and does not make use of closure under reversal.

\begin{figure}[t]
\begin{center}
\begin{picture}(400,120)
    \put(-30,0){\includegraphics{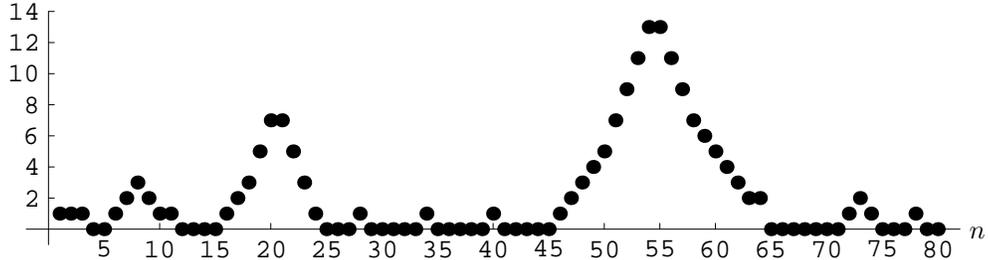}}
    \put(355,10){$n$}
\end{picture}
\end{center}
\caption{The multiplicity of 1 as an eigenvalue of $G_{2/(\sqrt 5-1)}(n)$\label{m(n)}}
\end{figure}

This result (and other computations) suggests that the eigenvalues of the Gram matrix
$G_\alpha(n):=(w_i \cdot w_j)_{1\le i,j\le n+1}$, where $F_n(W)=\{w_1,\dots,w_{n+1}\}$, may have
structure. Note that the eigenvalues of a Gram matrix do not depend on the ordering of the vectors,
and are necessarily nonnegative real numbers. A particularly striking phenomenon is the following.
Set $m(n)$ to be the multiplicity of 1 as an eigenvalue of $G_{2/(\sqrt{5}-1)}(n)$. For example
$m(55)=13$ and $m(65)=0$. Figure~\ref{m(n)} shows an impressive amount of structure, but this
author has no explanation for why {\em any} structure would exist as $n$ changes. Similar pictures
result from considering other irrationals.

\section{The route of the proof}
Let $\floor{x}$ denote the floor of $x$, and $\fp{x}$ the fractional part of $x$, i.e.,
$x=\floor{x}+\fp{x}$. Define $\B_\alpha(k):=\#\{q \colon 1\leq q < k,\,
\fp{q\alpha}<\fp{k\alpha}\}$, which counts the number of integers in $[1,k)$ that are `better'
denominators for approximating $\alpha$ from below.

Our proof proceeds by connecting the sum in Theorem~\ref{thm:main} to $\B_\alpha(n)$ (for some
$\alpha$), finding a recurrence satisfied by $\B_\alpha(n)$, and then reducing that recurrence
modulo 2.

The characteristic word with slope $\alpha$ is defined by
    \[c_\alpha(n):=\floor{(n+2)\alpha}-\floor{(n+1)\alpha}.\]
If $\alpha$ is irrational, then $c_\alpha$ is a Sturmian word~\cite{Lothaire}*{Thm 2.1.13}. It is
known \cite{Lothaire}*{Thm 2.1.3, Prop 2.1.18} that for every binary Sturmian word $W$ and natural
number $n$, there is an $\alpha\in(0,1)$ with $F_n(W)=F_n(c_\alpha)$, and so it suffices for our
purposes to consider characteristic words, and to write $F_n(\alpha):=F_n(c_\alpha)$.

\begin{lem}\label{lem:sourceofB}
\(\displaystyle \sum_{\vec w \in F_n(\alpha)} h(\vec w) = \B_\alpha(n)+(n+1)\floor{n\alpha}+1.\)
\end{lem}

\begin{lem}\label{Brepresentation}
Let $\alpha\in (0,1/2)$ be irrational. Then $\B_\alpha(k)+\B_{1-\alpha}(k)=k-1$. Moreover,
$\B_\alpha(1)=0$, $\B_\alpha(2)=1$, and for $k\ge3$
  $$
  \B_\alpha(k) -2\B_\alpha(k-1)+\B_\alpha(k-2) =
   \begin{cases}
    1-k,     & \fp{k\alpha}  \in [0,\alpha); \\
    k-1,     & \fp{k\alpha}  \in [\alpha,2\alpha); \\
    0,       & \fp{k\alpha}  \in [2\alpha,1).
   \end{cases}
  $$
\end{lem}

\begin{lem} \label{SignCorollary}
Let $\alpha\in(0,1)$ be irrational, and $k$ any positive integer. If $k$ is odd, then
$\B_\alpha(k)$ is even. If $k$ is even, then $\B_\alpha(k) \equiv \floor{k\alpha}+1 \pmod{2}$.
\end{lem}

\section{Proofs}
\begin{proof}[Proof of Lemma~\ref{lem:sourceofB}]
We begin by following \cite{Shallit}*{Lem 10.5.1}; define $\pi_i$ by
    \[\{0=\pi_0<\pi_1<\pi_2<\dots<\pi_n\}=\{0,\fp{-\alpha},\fp{-2\alpha},\dots,\fp{-n\alpha}\}.\]
Set $v_i(x):=\floor{(i+1)\alpha+x}-\floor{i \alpha+x}$, and set
    \[w_i:=\big(v_0(\pi_i),v_1(\pi_i),v_2(\pi_i),\dots,v_{n-1}(\pi_i)\big).\]
Nontrivially (see \cite{Shallit}),
    \(F_n(\alpha)=\big\{ w_i \colon
        0\le i \le n \big\}\),
and $w_0,w_1,\dots,w_n$ are ordered lexicographically. Elementary examination yields $h(w_i)=|\ZZ
\cap (\pi_i,n\alpha+\pi_i]|$, and this last quantity is either $\floor{n\alpha}$ or
$\floor{n\alpha}+1$. We start with $h(w_0)=\floor{n\alpha}$, and the first $i$ for which
$h(w_i)=\floor{n\alpha}+1$ is the $i$ for which $n\alpha+\pi_i\in\ZZ$, that is, when
$\fp{-n\alpha}=\pi_i$. In other words, the last $\B_\alpha(n)+1$ factors have weight
$\floor{n\alpha}+1$ and the first $n+1-(\B_\alpha(n)+1)$ factors have weight $\floor{n\alpha}$.
This gives
 \begin{multline*}
 \sum_{\vec{w} \in F_n(\alpha)} h(\vec{w}) =
                 (\B_\alpha(n)+1)(\floor{n\alpha}+1)
                 +(n+1-(\B_\alpha(n)+1))\floor{n\alpha}
            \\=  1+\B_\alpha(n)+(n+1)\floor{n\alpha}.
 \end{multline*}
\end{proof}

Our proof of Lemma~\ref{Brepresentation} is similar in spirit to, and was directly inspired by,
S\'{o}s's proof of the Three-Gap Theorem~\cite{Sos}.

\begin{proof}[Proof of Lemma~\ref{Brepresentation}]
Observe that $0<\fp{q\alpha} < \fp{k\alpha}$ iff $\fp{k(1-\alpha)}<\fp{q(1-\alpha)}<1$, so that $q$
with $1\leq q<k$ is in either the set $\{q: 1\leq q<k, \fp{q\alpha}<\fp{k\alpha}\}$ or in the set
$\{ q: 1\leq q<k, \fp{q(1-\alpha)}<\fp{k(1-\alpha)}\}$, and is not in both (as $\alpha$ is
irrational, $\fp{k(1-\alpha)}\not=\fp{q(1-\alpha)}$). Thus, $\B_\alpha(k)+\B_{1-\alpha}(k)=k-1$.

We think of the $k+2$ numbers $0,\fp{\alpha},\dots,\fp{k \alpha},1$ as lying on a unit circle, and
labeled $P_0, P_1, \dots, P_k, P_0$, respectively, i.e., $P_j := e^{2\pi j \alpha
\sqrt{-1}}=e^{2\pi \fp{j \alpha} \sqrt{-1}}$. ``The arc $\overline{P_iP_j}$'' refers to the
half-open counterclockwise arc from $P_i$ to $P_j$, containing $P_i$ but not $P_j$. We say that
three distinct points $A,B,C$ are {\em in order} if $B \not\in \overline{CA}$. We say that
$A,B,C,D$ are \emph{in order} if both $A,B,C$ and $C,D,A$ are in order. Essentially, if when moving
counter-clockwise around the circle starting from $A$, we encounter first the point $B$, then $C$,
then $D$, and finally $A$ (again), then $A,B,C,D$ are in order.

By rotating the circle through an angle of $2\pi\alpha$, so that $P_i \mapsto P_{i+1}$ ($0\leq i
\leq k$), we find that each $P$ on the arc $\overline{P_{k-2}P_{k-1}}$ is rotated onto a $P$ on the
arc $\overline{P_{k-1}P_k}$. Specifically, the number of $P_0,P_1,\dots,P_{k-2}$ on
$\overline{P_{k-2}P_{k-1}}$ is the same as the number of $P_1,P_2,\dots,P_{k-1}$ on
$\overline{P_{k-1}P_k}$. Set
 $$X:=\{P_0,P_1,\dots,P_{k-2}\} \quad\text{and}\quad
    Y:=\{P_1,P_2,\dots,P_{k-1}\},$$
so that what we have observed is
 \begin{equation}\label{dagger}
    \left| X \cap \overline{P_{k-2}P_{k-1}} \right| =
    \left| Y \cap \overline{P_{k-1}P_k} \right|.
 \end{equation}
Also, we will use the definition of $\B_\alpha$ in the forms
 $ \B_\alpha(k) = \left| Y \cap \overline{P_0P_k}\right|$
and
    $\B_\alpha(k-1)=|X \cap \overline{P_0P_{k-1}}| -1$,
and with $k$ and $k-1$ replaced by $k-1$ and $k-2$, when circumstances allow.

Now, first, suppose that $\fp{k\alpha} \in [0,\alpha)$, so that the points $P_0, P_k, P_{k-2},
P_{k-1}$ are in order on the circle. We have
  \begin{align*}
    X \cap \overline{P_{k-2}P_{k-1}}
        &=  X \cap \left( \overline{P_0P_{k-1}} \setminus
                \overline{P_0P_{k-2}}\right)\\
        &=  \left( X \cap \overline{P_0P_{k-1}}\right) \setminus
                \left(X \cap \overline{P_0P_{k-2}}\right) \\
    \left|X \cap \overline{P_{k-2}P_{k-1}}\right|
        &=  \left|\left( X \cap \overline{P_0P_{k-1}}\right) \right|\,-\,
                \left| \left(X \cap \overline{P_0P_{k-2}}\right)\right| \\
        &=  \left(\B_\alpha(k-1)-1\right)-\left(\B_\alpha(k-2)-1\right)\\
        &=  \B_\alpha(k-1)-\B_\alpha(k-2),
  \end{align*}
and similarly
  \begin{align*}
    Y \cap \overline{P_{k-1}P_k}
        &=  \left( Y \cap \overline{P_{k-1}P_0}\right) \cup
                \left( Y \cap \overline{P_0P_k}\right)\\
        &=  \left( Y \setminus \left( Y \cap \overline{P_0P_{k-1}} \right) \right)
                \cup \left( Y \cap \overline{P_0P_k} \right) \\
    \left| Y \cap \overline{P_{k-1}P_k} \right|
        &=  (|Y|-\left| Y \cap \overline{P_0P_{k-1}} \right|)\,+\,
                \left|Y \cap \overline{P_0P_k} \right| \\
        &=  (k-1-\B_\alpha(k-1))+\B_\alpha(k)
  \end{align*}
so that Eq.~\eqref{dagger} becomes $\B_\alpha(k-1)-\B_\alpha(k-2)=\B_\alpha(k)-\B_\alpha(k-1)+k-1$,
as claimed in the statement of this lemma.

Now suppose that $\fp{k\alpha} \in [\alpha,2\alpha)$, so that the points $P_0, P_{k-1}, P_k,
P_{k-2}$ are in order. By arguing as in the above case, we find
    $$X \cap \overline{P_{k-2}P_{k-1}}=\left(X\setminus\left(X\cap
    \overline{P_0P_{k-2}}\right)\right)\cup \left(X\cap \overline{P_0P_{k-1}}\right),$$
and so $\left|X\cap\overline{P_{k-2}P_{k-1}}\right|=k-1-(\B(k-2)-1)+(\B(k-1)-1)$. Likewise,
    $$Y \cap \overline{P_{k-1}P_k}=\left(Y\cap\overline{P_0P_k}\right)\setminus
    \left(Y\cap\overline{P_0P_{k-1}}\right)$$
so that $\left|Y \cap \overline{P_{k-1}P_k} \right|=\B_\alpha(k)-\B_\alpha(k-1)$. Thus, in this
case Eq.~\eqref{dagger} becomes $\B_\alpha(k-1)-\B_\alpha(k-2)+k-1=\B_\alpha(k)-\B_\alpha(k-1)$, as
claimed in the statement of the lemma.

Finally, suppose that $\fp{k\alpha} \in [2\alpha,1)$, so that the points $P_0,P_{k-2},P_{k-1},P_k$
are in order. We find
 $$X\cap\overline{P_{k-2}P_{k-1}}=\left( X \cap \overline{P_0P_{k-1}}\right)
 \setminus \left(X \cap \overline{P_0P_{k-2}}\right)$$
and so $\left|X\cap\overline{P_{k-2}P_{k-1}}\right|=\B(k-1)-\B(k-2)$. Also,
 $$Y\cap\overline{P_{k-1}P_k} = \left(Y\cap\overline{P_0P_k}\right)\setminus
 \left(Y\cap\overline{P_0P_{k-1}}\right)$$
and so $\left|Y\cap\overline{P_{k-1}P_k}\right|=\B_\alpha(k)-\B_\alpha(k-1)$. As claimed,
Eq.~\eqref{dagger} becomes $\B_\alpha(k-1)-\B_\alpha(k-2)=\B_\alpha(k)-\B_\alpha(k-1)$.
\end{proof}

For the remaining proofs, we write
 $\displaystyle \tf{Q}:=\left\{%
\begin{array}{ll}
    1, & \hbox{$Q$ is true;} \\
    0, & \hbox{$Q$ is false.}
\end{array}%
\right.$

\begin{proof}[Proof of Lemma~\ref{SignCorollary}]
Reducing Lemma~\ref{Brepresentation} modulo 2, we find that if $0<\alpha<1/2$, then
    \begin{equation*}\label{Bmod2}
    \B_\alpha(k) \equiv \B_\alpha(k-2) +\tf{k \text{ even}}\,\tf{\fp{k\alpha}<2\alpha}.
    \end{equation*}
and if $1/2<\alpha<1$, then
    \begin{equation*}
    B_\alpha(k)=-B_{1-\alpha}(k)+k-1\equiv B_{1-\alpha}(k)+\tf{k\text{ even}} \pmod 2
    \end{equation*}
We work in four cases: $k$ may be odd or even, and $\alpha$ may be less than or greater than $1/2$.

Assume first that $k$ is odd and $0<\alpha<1/2$. As $\B_\alpha(1)=0$ and $\B_\alpha(k)\equiv
\B_\alpha(k-2)\pmod2$, we see by induction that $\B_\alpha(k)$ is even.

Now assume that $k$ is odd and $1/2<\alpha<1$. We have $\B_\alpha(k)\equiv \B_{1-\alpha}(k)\pmod2$,
and as $0<1-\alpha<1/2$, the paragraph immediately above implies that $\B_{1-\alpha}(k)$ is even.

Now assume that $k$ is even and $0<\alpha<\tfrac12$. Set $\beta=2\alpha$, $k=2\ell$ and
$\BB(i)=\B_\alpha(2i)$, so that $\BB(1)=1$ and
    $\BB(i)\equiv \BB(i-1)+\tf{\fp{i\beta} <\beta} \pmod2.$
We have
 \begin{align*}
    \B_\alpha(k)=\B_\alpha(2\ell)
                =      \BB(\ell)
                &\equiv \BB(\ell-1)+\tf{\fp{\ell\beta}<\beta}
                                \pmod{2}\\
                &\equiv \BB(1)+\sum_{i=2}^{\ell}
                         \tf{\fp{i\beta} < \beta} \pmod{2}\\
                &=      \BB(1)+\floor{\ell\beta}
                =      1+\floor{(k/2)(2\alpha)} = 1+\floor{k\alpha},
 \end{align*}
since $\sum_{i=2}^{\ell} \tf{\fp{i\beta} <\beta}$, with $\beta\in(0,1)$, counts the integers in the
interval $(\beta,\ell\beta]$.

Finally, suppose that $k$ is even and $1/2<\alpha<1$. By the paragraph immediately above,
$\B_{1-\alpha}(k) \equiv 1+\floor{k(1-\alpha)}\pmod{2}$. We have
 \begin{align*}
  \B_\alpha(k)
        &\equiv 1+B_{1-\alpha}(k) \pmod{2} \\
        &\equiv \floor{k(1-\alpha)}\pmod{2} \\
        &\equiv \floor{k\alpha}+1 \pmod{2},
 \end{align*}
where we have again used the irrationality of $\alpha$ in the last line.
\end{proof}

\begin{proof} [Proof of Theorem \ref{thm:main}]
Now, if $n$ is odd, then by Lemma~\ref{SignCorollary}, $\B_\alpha(n)\equiv0\pmod{2}$ and obviously
$n+1\equiv 0 \pmod{2}$, whence $1+\B_\alpha(n)+(n+1)\floor{n\alpha} \equiv 1 \pmod{2}$. If $n$ is
even, then by Lemma~\ref{SignCorollary}, $\B_\alpha(n)\equiv \floor{n\alpha}+1$, whence
$1+\B_\alpha(n)+(n+1)\floor{n\alpha} \equiv 0 \pmod{2}$.
\end{proof}

\end{document}